\documentclass{ifacconf}

\makeatletter
\let\old@ssect\@ssect %
\makeatother

\usepackage{natbib}        %

\usepackage[short]{optidef}
\usepackage{tikz, lscape, amsmath,amsfonts,latexsym,amssymb}
\usepackage[utf8]{inputenc}
\usepackage[english]{babel}
\usepackage{lmodern, verbatim}
\usepackage[T1]{fontenc}
\usepackage{calc,bm}
\usepackage{enumitem,nameref,hyperref}
\usepackage[mathscr]{eucal}
\usepackage{stackrel,floatrow}

\makeatletter
\def\@ssect#1#2#3#4#5#6{%
	\NR@gettitle{#6}%
	\old@ssect{#1}{#2}{#3}{#4}{#5}{#6}%
}
\makeatother

\renewcommand{\b}{\mathbf}

\newcommand{\barRho}{\bar{\rho}_\epsilon}
\newcommand{\barX}{\bar{\b x}}
\newcommand{\barU}{\bar{\b u}}

\newcommand{\BB}{\mathbb{B}}

\newcommand{\epsArr}{\bm{\epsilon}}%

\newcommand{\NN}{\mathbb{N}}

\newcommand{\RR}{\mathbb{R}}

\newcommand{\A}{\mathcal{A}}

\newcommand{\C}{\mathcal{C}}

\DeclareMathOperator*{\argmin}{arg\,min}

\DeclareMathOperator*{\Int}{Int}

\newcommand{\tb}{\textbf} %
\newcommand{\R}{\RR} %
\renewcommand{\b}{\mathbf}

\newcommand{\iv}[2]{[\![#1,#2]\!]} %

\date{\today}

\makeatletter
\let\orgdescriptionlabel\descriptionlabel
\renewcommand*{\descriptionlabel}[1]{%
	\let\orglabel\label
	\let\label\@gobble
	\phantomsection
	\edef\@currentlabel{#1\unskip}%
	\let\label\orglabel
	\orgdescriptionlabel{#1}%
}
\makeatother

\begin{document}
	\begin{frontmatter}
		
		\title{Stability of solutions for controlled nonlinear systems under perturbation of state constraints}

		\author[First]{Pierre-Cyril Aubin-Frankowski}
		
		\address[First]{INRIA Paris, France (e-mail: pierre-cyril.aubin@inria.fr).}
		
		\begin{abstract}: This paper tackles the problem of nonlinear systems, with sublinear growth but unbounded control, under perturbation of some time-varying state constraints. It is shown that, given a trajectory to be approximated, one can find a neighboring one that lies in the interior of the constraints, and which can be made arbitrarily close to the reference trajectory both in $L^\infty$-distance and $L^2$-control cost. This result is an important tool to prove the convergence of approximation schemes of state constraints based on interior solutions and is applicable to control-affine systems.
		\end{abstract}
		
		
	\end{frontmatter}
	\section{Introduction}
	We consider a nonlinear system with unbounded control and state constraints
	\begin{align}
	\b x'(t)&=\b f(t,\b x(t),\b u(t)), &&\text{ for a.e.\ } t\in[0,T], \label{eq_def_system} \\
	\b x(t)&\in\A_{0,t}:=\{\b x\,|\, \b h(t,\b x)\le 0\},&&\text{ for all } t\in[0,T], \label{eq_def_constraints}
	\end{align}
	where $\b f:[0,T]\times\R^{N}\times\R^{M}\rightarrow \R^N$ and $\b h:[0,T]\times\R^{N}\rightarrow \R^P$. Given a reference trajectory $\barX(\cdot)$, such that $\b h(0,\barX(0))<0$, with control $\barU(\cdot)$ satisfying \eqref{eq_def_system}-\eqref{eq_def_constraints}, our goal is to design a trajectory $\b x^\epsilon(\cdot)$ with the same initial condition and some control $\b u^\epsilon(\cdot)$, chosen such that $\b x^\epsilon(\cdot)$ can be made arbitrarily $L^\infty$-close to $\barX(\cdot)$, with $\b u^\epsilon(\cdot)$ having almost the same $L^2$-norm as $\barU(\cdot)$, while also satisfying \eqref{eq_def_system} and the following tightened constraints:%
	\begin{equation}\label{eq_def_pertubed_constraints}
	\b x^\epsilon(t)\in\A_{\epsilon,t}:=\{\b x\,|\, \epsArr+\b h(t,\b x)\le 0\}\text{ for all } t\in[0,T].
	\end{equation}
	This construction is crucial to prove the convergence of approximation schemes of the constraints from within, in the sense that if $(\barX(\cdot),\barU(\cdot))$ is the solution of some optimal control problem with quadratic cost in control, then $(\b x^\epsilon(\cdot),\b u^\epsilon(\cdot))$ would be almost optimal while strictly interior. Such schemes were used by the author in \cite{aubin2020hard_control} for linear $\b f$ and $\b h$, leveraging convexity of the set of trajectories for which \eqref{eq_def_system}-\eqref{eq_def_constraints} hold. Here we provide instead assumptions on $\b f$ and $\b h$ designed originally by \cite{Bettiol2012LEO} for bounded differential inclusions with time-invariant constraints. We further improve on their construction to have both an estimate on the $L^2$-norm and to cover unbounded systems \eqref{eq_def_system} and time-varying constraints \eqref{eq_def_constraints}. This analysis can be related also to \cite{Bettiol2011} where time-dependent bounded systems are considered. The prototypical cases we are interested in are constrained nonlinear control-affine systems as studied e.g.\ in a more restrictive setting in \cite[Section 4]{Cannarsa2008} or \cite{Aronna2016}.	
	
	\section{Main result}

	\tb{Notations.} The integer interval is written $\iv{i}{j}=\{i,i+1,\dots,j\}$. We denote by $\R_+$ the set of nonnegative reals, and use the shorthand $L^p(0,T)$ for $L^p([0,T],\R^d)$ with $p\in\{1,2,\infty\}$, and $L^p_+$ when the output set is $\R_+^d$. The set $\BB_d$ is the closed Euclidean unit ball of $\R^d$ of center $\b 0$. Given a set $\A\subset \R^d$, $\Int(\A)$ designates its interior, $\partial \A$ its boundary, and $d_\A(\cdot)$ is the Euclidean distance to $\A$.  
	
	\begin{figure}
		\begin{center}
			\includegraphics[width=8.4cm]{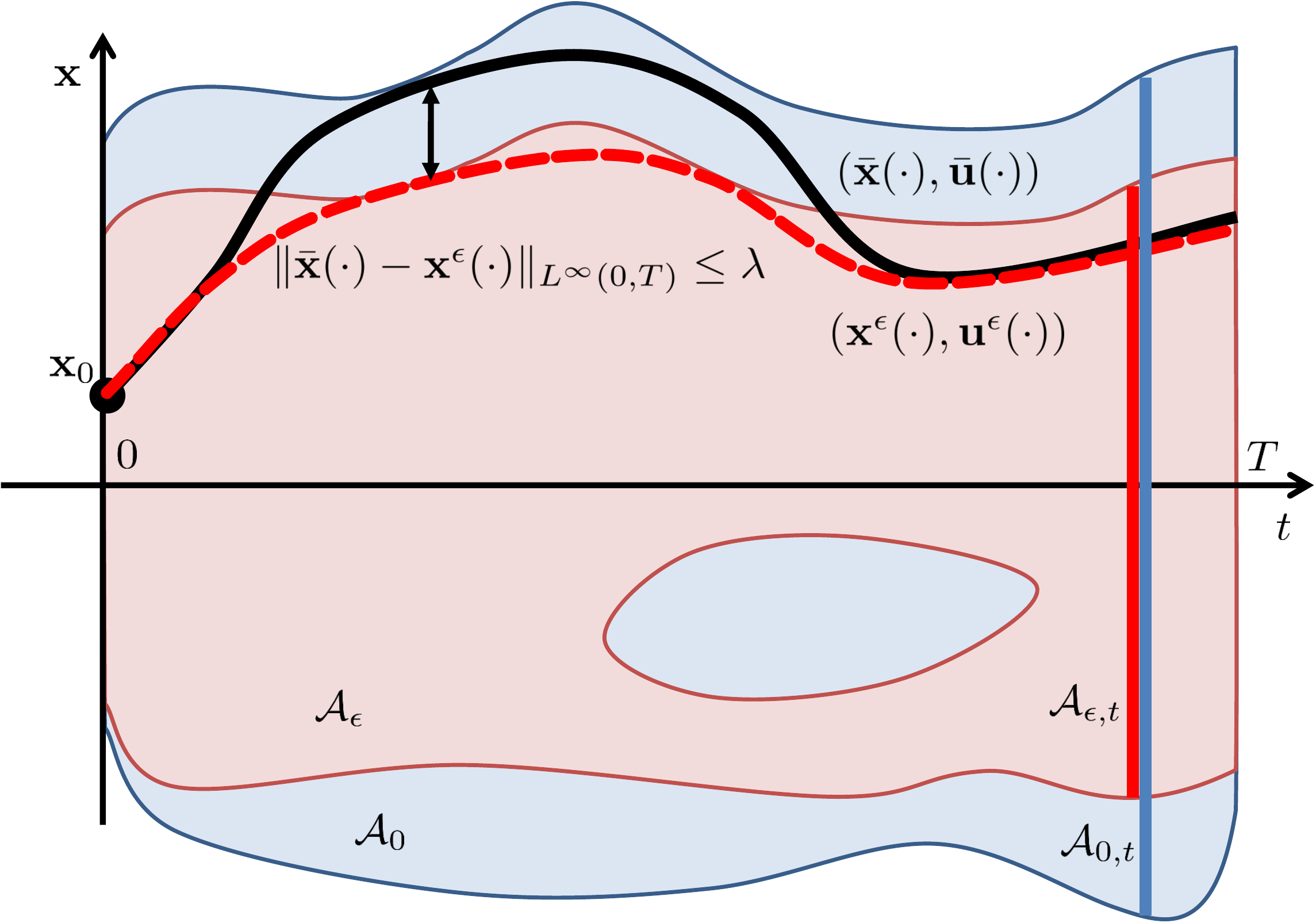}    %
			\caption{Illustration of the trajectories and constraints considered in Theorem \ref{thm_NFT_unbounded}.} 
			\label{fig:schema}
		\end{center}
	\end{figure}

	We call $\b f$-trajectories the solutions of \eqref{eq_def_system} for measurable controls $\b u(\cdot)$. For any $\epsArr\in\R_+^P$, define  $$\A_\epsilon:=\{(t,\b x)\,|\,t\in[0,T],\, \b x\in\A_{\epsilon,t}\}.$$ 
	A trajectory is said to be $\A_\epsilon$-feasible if \eqref{eq_def_pertubed_constraints} holds, for instance $\barX(\cdot)$ is $\A_0$-feasible by assumption. We define the maximal constraint violation $\rho_{\epsilon,[t_0,t_1]}(\b x(\cdot))$ of a trajectory on an interval $[t_0,t_1]\subset [0,T]$ as follows%
	\begin{gather*}
	\rho_{\epsilon,[t_0,t_1]}(\b x(\cdot)):=\sup _{t \in [t_0,t_1]} d_{\A_{\epsilon,t}}(\b x(t)).
	\end{gather*}
	We assume from now on that the control $\barU(\cdot)$ belongs to $L^\infty(0,T)$. This restriction is due to \ref{hyp_lipF_nonlin} below since we use a time-delay in the construction of the control $\b u^\epsilon(\cdot)$ that is ill-suited to track the distance between controls. If $\b f(t,\cdot,\b u)$ is $k_f(t)$-Lipschitz for any $\b u\in\R^M$, then we may just assume that $\barU(\cdot)\in L^2(0,T)$.

	\begin{description}
		\item [(H-1)\label{hyp_regPerturbation_nonlin}] (Regular perturbation of $\A$) 
		\begin{gather*}
		\hspace{-5mm}\forall\, \lambda>0,\;\exists\,  \epsArr>0,\\
		 \forall\,(t,\b x)\in \A_0\cap([0,T]\times \|\barX(\cdot)\|_{L^{\infty}(0,T)}\BB_N),\\ d_{\A_{\epsilon,t}}(\b x)\leq \lambda.
		\end{gather*}
		\item  [(H-2)\label{hyp_C0_dA_nonlin}] (Uniform continuity from the right of $d_{\partial\A_{\epsilon,t}}$ w.r.t. $\epsilon$ and $t$)
		There exist $\epsArr_0>0$, $\Delta_0>0$, and a continuous function $\omega_\A(\cdot)\in\C^0(\R_+,\R_+)$ such that $\omega_\A(0)=0$ and, for all $\epsArr\leq\epsArr_0$,  and all $(t,\b x)\in \A_0\cap([0,T]\times 2\|\barX(\cdot)\|_{L^{\infty}(0,T)}\BB_N)$,
		\begin{equation*}
		\hspace{-7mm}\forall\, \delta \in[0,\min(\Delta_0,T-t)],\; \|d_{\partial\A_{\epsilon,t+\delta}}(\b x)-d_{\partial\A_{\epsilon,t}}(\b x)\|\leq \omega_\A(\delta).
		\end{equation*}\\
		\item  [(H-3)\label{hyp_subLin_nonlin}] (Sublinear growth of $\b f$ w.r.t. $\b x$ and $\b u$)
		\begin{gather*}
		\exists\, \theta(\cdot)\in L_+^2(0,T),\; \forall\, t\in[0,T],\; \forall \, \b x\in \R^N,\; \forall \, \b u\in \R^M,\\
		\|\b f(t,\b x,\b u)\|\leq \theta(t)(1+\|\b x\|+\|\b u\|).
		\end{gather*}\\
		\item  [(H-4)\label{hyp_IP_Keps_nonlin}] (Inward-pointing condition) There exist $\epsArr_0>0$, $M_u>0$, $M_v>0$, $\xi>0$, and $\eta>0$ such that for all $\epsArr\le\epsArr_0$ and all $(t,\b x)\in (\partial \A_\epsilon +(0,\eta\BB_N)) \cap \A_\epsilon\cap ([0,T]\times (1+2\|\barX(\cdot)\|_{L^{\infty}(0,T)})\BB_N)$, we can find $\b u\in M_u\BB_M$ such that $\b v:=\b f(t,\b x,\b u)$ belongs to $M_v \BB_N$ and 
		\begin{equation}\label{eq_IPC_nonlin}
		\b y+\delta(\b v+\xi\BB_N)\subset\A_{\epsilon, t+\delta}
		\end{equation}
		for all $\delta\in[0,\xi]$ and all $\b y\in (\b x+\xi\BB_N)\cap\A_{\epsilon, t}$.\\ %
		\item  [(H-5)\label{hyp_absC_nonlin}] (Left local absolute continuity of $\b f$ w.r.t. $t$)
		\begin{gather*}%
		\exists\, \gamma(\cdot) \in L_+^{1}(0,T),\; \exists\, \beta_u(\cdot) \in L_+^{2}(0,T),\\
		 \forall\,0\le s < t\le T ,\; \forall \, \b x\in(1+2\|\barX(\cdot)\|_{L^{\infty}(0,T)})\BB_N,\nonumber\\
		\hspace{-.5cm}\forall \, \b u_s\in (M_u+\|\barU(s)\|)\BB_M,\; \exists \, \b u_t\in \b u_s + \beta_u(s)\BB_M,\\
		\|\b f(t,\b x,\b u_t)-\b f(s,\b x,\b u_s)\|\leq \int_{s}^{t} \gamma\left(\sigma\right) d \sigma.%
		\end{gather*}
		
		\begin{align}
		\hspace{-1.2cm}\text{Let } &R:=e^{\|\theta(\cdot)\|_{L^{1}(0,T)}}\left[1+\|\barX(\cdot)\|_{L^{\infty}(0,T)}\right.\nonumber\\
		&\hspace{-5mm}+(1+M_u)\|\theta(\cdot)\|_{L^{1}(0,T)}\nonumber\\
		&\hspace{-5mm}+
		\left.\|\theta(\cdot)\|_{L^{2}(0,T)}(\|\barU(\cdot)\|_{L^{2}(0,T)}+\|\beta_u(\cdot)\|_{L^{2}(0,T)})\right] . \label{eq_def_R}
		\end{align}\\
		\item  [(H-6)\label{hyp_lipF_nonlin}] (Local Lipschitz continuity of $\b f$ w.r.t. $\b x$)
		\begin{gather*}
		\exists\, k_f(\cdot)\in L_+^2(0,T),\; \forall\, t\in[0,T],\; \forall \, \b x,\b y\in R\BB_N,\\ \forall \, \b u\in (M_u+\|\barU(\cdot)\|_{L^\infty(0,T)})\BB_M,\\
		\|\b f(t,\b x,\b u)-\b f(t,\b y,\b u)\|\leq k_f(t)\|\b x-\b y\|.
		\end{gather*}\\
		\item  [(H-7)\label{hyp_lipSelection_nonlin}] (Hölderian selection of the controls in \ref{hyp_absC_nonlin})
		\begin{gather*}
		\exists\, \gamma(\cdot) \in L_+^{1}(0,T),\; \exists\, \alpha\in]0,1],\; \exists\, k_u(\cdot) \in L_+^{2}(0,T),\\ \forall\, 0\le s < t\le T,\;
		\forall \, \b x\in(1+2\|\barX(\cdot)\|_{L^{\infty}(0,T)})\BB_N,\\
		\hspace{-5mm}\forall \, \b u_s\in (M_u+\|\barU(s)\|)\BB_M,\;
		\exists \, \b u_t\in \b u_s+(t-s)^\alpha k_u(s)\BB_M,\\
		\|\b f(t,\b x,\b u_t)-\b f(s,\b x,\b u_s)\|\leq \int_{s}^{t} \gamma\left(\sigma\right) d \sigma.%
		\end{gather*}
	\end{description}
	
	\begin{thm}\label{thm_NFT_unbounded} Under assumptions \ref{hyp_regPerturbation_nonlin}-\ref{hyp_lipF_nonlin}, for any $\lambda>0$, there exists $\epsArr>0$ and a $\b f$-trajectory $\b x^\epsilon(\cdot)$ on $[0,T]$ such that $\b x^\epsilon(0)=\barX(0)$, $\b x^\epsilon(t) \in \Int \A_{\epsilon,t}$ for all $t \in[0,T]$, and $$\|\barX(\cdot)-\b x^\epsilon(\cdot)\|_{L^{\infty}\left(0,T\right)} \leq \lambda.$$ 
		Moreover if \ref{hyp_lipSelection_nonlin} is satisfied, then, for any mapping $\b R(\cdot)\in\C^0([0,T],\R^{M,M})$ with positive semidefinite matrix values, one can choose $\epsArr>0$ and $\b x^\epsilon(\cdot)$ such that the controls $\b u^\epsilon(\cdot)$ satisfy 
		$$\left\vert\Vert \b R(\cdot)^{1/2} \barU(\cdot)\Vert^2_{L^2(0,T)}-\Vert \b R(\cdot)^{1/2}  \b u^\epsilon(\cdot)\Vert^2_{L^2(0,T)}\right\vert\leq\lambda.$$
	\end{thm}

	\noindent\textbf{Discussion of the assumptions:} The properties \ref{hyp_regPerturbation_nonlin} and \ref{hyp_C0_dA_nonlin} imposed on the constraint set can for instance be derived from the $\C^{1,1}$-regularity of $\b h$, coupled with the assumptions that the Jacobian $\frac{\partial \b h(t,\b  x)}{\partial\b  x}$ of $\b h$ at all $(t,\b x)\in \partial \A_0$ is surjective. The classical growth assumption \ref{hyp_subLin_nonlin} prevents in turn finite-time explosion of the trajectories. The Lipschitzianity \ref{hyp_lipF_nonlin} guarantees their uniqueness and was designed to encompass control-affine systems of the form $\b x'(t)=\b a(t,\b x)+\b b(t,\b x)\b u$ with $\tilde{k}_f(t)$-Lipschitz functions $\b a(t,\cdot)$ and $\b b(t,\cdot)$, for some $\tilde{k}_f(\cdot)\in L^2(0,T)$. The other assumptions are more technical and inspired by \cite{Bettiol2012LEO}. Inward-pointing conditions such as \ref{hyp_IP_Keps_nonlin}, which can be deduced from a normal cone formulation \cite[Lemma 5.3]{Bettiol2012LEO}, have been shown to yield the $L^\infty$-bounds we seek. The time regularity \ref{hyp_absC_nonlin} was introduced to tackle discontinuities in the dynamics, and showcased on a civil engineering example \cite[Section 4]{Bettiol2012LEO}. We adapt it to control systems and refine it in \ref{hyp_lipSelection_nonlin}. 
	
	Note that state constraints of order 2 (or more), e.g.\ $\ddot{x}=u$ with $x$ constrained, do not enter into the proposed framework as the inward-pointing assumption does not hold in these cases, being limited to ``order 1 constraints''. The proof can in principle be adapted to systems $\tilde{\b f}$ with control constraints following a Lipschitz (or Hölderian) closed-valued map $t\leadsto U(t)$ by considering the projection over $U(t)$ and $\b f(t,\b x,\b u)=\tilde{\b f}(t,\b x,\operatorname{proj}_{U(t)}(\b u))$ assuming that this $\b f$ satisfies the above assumptions.
	
	\noindent\textbf{Idea of the proof:}  The overall strategy to construct a neighboring $\A_\epsilon$-feasible trajectory can be related to that of \cite{Bettiol2012LEO}. Modifying it to unbounded controls and time-varying constraints is however not straightforward. We start by considering small subintervals $[0,T]=\bigcup_{i\in[\![0,N_0-1]\!]} [t_i,t_{i+1}]$ and proceed iteratively. If the $i$th-trajectory stays in $\A_\epsilon$ over $[t_i,t_{i+1}]$, we move to the next time interval. Otherwise for the $(i+1)$th-trajectory over $[t_i,t_{i+1}]$, \ref{hyp_IP_Keps_nonlin} provides us with an inward-pointing control to stay in $\A_\epsilon$ for a short time. Then we apply a delayed control of the $i$th-trajectory for the rest of $[t_i,t_{i+1}]$, and the original control $\barU(\cdot)$ over $[t_{i+1},T]$. By monitoring several quantities, we can show that the resulting control after $N_0$ iterations is $L^2$-close from  $\barU(\cdot)$ and that the obtained trajectory is in $\A_\epsilon$.\\
	
	\noindent\textbf{Example:} Consider an electric motor
	\begin{equation*}
	x'(t)=a(t,x)+b(t,u),
	\end{equation*} with a bounded $a\in C^{1,1}([0,2]\times \R,\R)$ and constraints $h(x)=1-|x|$, for controls $u\in\R$. The motor suffers an incident at $T=1$. If it is a power surge $$b(t,u)=\tilde{b}(t)u=\left\{
	\begin{array}{ll}
	u &\text{ if $t\in[0,1]$}\\
	u/\sqrt[4]{t-1} &\text{ if $t\in]1,2]$}
	\end{array}
	\right.,$$ then \ref{hyp_subLin_nonlin} holds for $\theta\equiv \tilde{b}+\|f\|_\infty$. It remains to check \ref{hyp_lipSelection_nonlin}. For $s < t\le 1 $, take $u_t=u_s$. For $s< 1 <  t $, set $u_t=\sqrt[4]{t-1}u_s$, so $|u_t-u_s|\le |u_s| \le \frac{\sqrt[4]{t-s}}{\sqrt[4]{|1-s|}}|u_s|$. For $s= 1 <  t $, take also $u_t=\sqrt[4]{t-1}u_s$, thus $|u_t-u_s|\le |u_s|$. For $1< s\le t $, set $u_t=\frac{\sqrt[4]{t-1}}{\sqrt[4]{s-1}}u_s$, so, by subadditivity of $\sqrt[4]{\cdot}$, $$|u_t-u_s| \le \frac{\sqrt[4]{t-1}-\sqrt[4]{s-1}}{\sqrt[4]{|1-s|}}|u_s|\le \frac{\sqrt[4]{t-s}}{\sqrt[4]{|1-s|}}|u_s|$$. Hence  \ref{hyp_lipSelection_nonlin} is indeed satisfied for $\gamma \equiv 0$, $k_u(s) =\frac{1}{\sqrt[4]{|1-s|}}(M_u+|\barU(s)|)$ and $\alpha=\frac{1}{4}$. If the incident consists in a power decline 
	$$b(t,u)=\left\{
	\begin{array}{ll}
	\arctan(u) &\text{ if $t\in[0,1]$}\\
	(1-\frac{\sqrt{t-1}}{2})\arctan(u) &\text{ if $t\in]1,2]$}
	\end{array}
	\right.,$$
	 then the system is bounded and \ref{hyp_lipSelection_nonlin} holds with $u_t=u_s$, $\gamma(\sigma)=\frac{1}{4\sqrt{\sigma-1}}$ for $\sigma\in]1,2]$ and $\gamma(\sigma)=0$ otherwise. In both cases \ref{hyp_IP_Keps_nonlin} is satisfied, so perturbing the constraints still allows for a trajectory and control close to the reference ones as per Theorem \ref{thm_NFT_unbounded}.\\

	\section{Proof of Theorem 1}
	 Below when we write $[s,t]\subset[0,T]$, we mean $0\le s<t\le T$. Let $\gamma(\cdot)$, $\beta_u(\cdot)$, $M_u$, $M_v$, $\xi$, $\eta$, $\Delta_0$ and $\omega_\A(\cdot)$ be as in \ref{hyp_C0_dA_nonlin}, \ref{hyp_IP_Keps_nonlin},  and \ref{hyp_absC_nonlin}. Define the following three non-negative continuous functions for $\delta \in[0,T]$ with limit $0$ at $\delta=0^+$:
	\begin{gather}
	\omega_\gamma(\delta):=\sup_{\substack{[s,t]\subset[0,T]\\|s-t|\leq \delta}} \left\{\int_s^t \gamma(\sigma) d \sigma\right\} \,; \label{eq_def_omega} \\ \omega_f(\delta):=\sup_{\substack{[s,t]\subset[0,T]\\|s-t|\leq \delta}}\left\{\int_s^t k_f(\sigma) d \sigma\right\} \nonumber \\
	\bar{\omega}(\delta):=\sup_{\substack{[s,t]\subset[0,T]\\|s-t|\leq \delta}} \left[\|\barX(t)-\barX(s)\|+(R+M_u) \|\theta(\cdot)\|_{L^1(s,t)}\right. \nonumber 
	\\	 \hspace{1cm}
	\left.+\|\theta(\cdot)\|_{L^2(s,t)} \left(\|\barU(\cdot)\|_{L^{2}(0,T)}+\|\beta_u(\cdot)\|_{L^{2}(0,T)}\right)\right]. \nonumber
	\end{gather}\\

	By definition of $R$ in \eqref{eq_def_R}, we have that $\barX([t_0,T])\subset (R-1)\BB_N$. Let $\Delta>0$, $k>0$ and $\hat{\rho}>0$, chosen in this order, be such that 
	\begin{gather}
	\Delta\leq \min(\xi,\Delta_0) \,;\, \omega_\A(\Delta)\leq \eta/4 \,;\, \bar{\omega}(\Delta)\leq \eta/4 \label{eq_def_Delta} \\
	M_{\Delta}:=\omega_\gamma(\Delta) + M_v\omega_f(\Delta) \,;\, C_{v,\Delta}:=M_v+M_{\Delta} e^{\omega_f(\Delta)}\label{eq_def_csts_nonlin} \\
	k\hat{\rho}M_v\le 1 \,;\, M_{\Delta} e^{\omega_f(\Delta)}\le \frac{\xi}{2} \,;\, 2<k\xi \,; \label{eq_hyp_csts_nonlin}\\
	 1 + k M_\Delta (1+e^{\omega_f(\Delta)}) e^{\omega_f(\Delta)}\leq \frac{k\xi}{2}.\nonumber
	\end{gather}\\

	One can easily check that these assumptions are compatible. The variables $\Delta$ and $k\hat{\rho}$ will serve as durations, whereas $\hat{\rho}$ is a bound on distances between trajectories. Let $N_0$ be the smallest integer such that $T\leq N_0 \Delta$. Fix $(t_i)_{i\in[\![0,N_0]\!]}$ with $0< t_{i+1}-t_i\leq \Delta$ such that 
	$$ [0,T]=\bigcup_{i\in[\![0,N_0-1]\!]} [t_i,t_{i+1}].$$ 
	Let
	\begin{align*}
		g(\rho)&:=e^{\omega_f(T)} \left[\bar{\omega}(k\rho)+ k\rho (C_{v,\Delta}+ M_\Delta e^{2\omega_f(\Delta)})\right];\\
		\tilde{g}(\rho)&:=\rho + g(\rho).
	\end{align*}
	For $g:\R_+\rightarrow \R_+$, we use the notation $g^{\circ i}(\cdot)$ for $g$ composed $i$ times with itself.
	Define, for any $n\in\NN^*$, $\tilde{d}_{n}(\rho):=\sum_{i=1}^{n}  \tilde{g}^{\circ i}(\rho)$. The functions $g$, $\tilde{g}$, $\tilde{d}_{n}$ are monotonically increasing, and $(\tilde{g}^{\circ n}(\cdot))_{n\in\NN}$ and $(\tilde{d}_{n}(\cdot))_{n\in\NN}$ form increasing sequences. By \ref{hyp_regPerturbation_nonlin}, there exists $\epsArr$ small enough as to have both $\A_{\epsilon,t}\neq \emptyset$ for all $t$, and $\barX(0) \in \Int \A_{\epsilon,0}$, while also satisfying, for $\barRho:=\rho_{\epsilon,[0,T]}(\barX(\cdot))$,
	\begin{equation} \label{eq_def_dN}
	\tilde{d}_{N_0}(\barRho)\leq \min(\lambda,\|\barX(\cdot)\|_{L^{\infty}(0,T)}) \,;\, \tilde{g}^{\circ N_0}(\barRho)\leq \hat{\rho} \,;\, \barRho\le \hat{\rho} .
	\end{equation}\\
	
	\noindent\textbf{Step 1}: starting from $t_0=0$, we construct a $\b f$-trajectory $\b x_1(\cdot)$ such that, for all $t\in[t_0,t_1]$, $\b x_1(t)\in \Int\A_{\epsilon,t}$. We distinguish two cases, whether $\barX(t_0)$ is far or close to the boundary $\partial\A_{\epsilon,t_0}$.\\
	
	\noindent\textit{Case 1}: If $d_{\partial\A_{\epsilon,t_0}}(\barX(t_0))>\eta/2$, then for any $t\in[t_0,t_1]$, by definition of $\omega_\A(\cdot)$ and $\bar{\omega}(\cdot)$,
	\begin{align*}
	d_{\partial\A_{\epsilon,t}}(\barX(t))&\geq d_{\partial\A_{\epsilon,t}}(\barX(t_0))-\|\barX(t)-\barX(t_0)\|\\
	&\geq d_{\partial\A_{\epsilon,t_0}}(\barX(t_0))- \omega_\A(\Delta)-\bar{\omega}(\Delta)\stackrel{\eqref{eq_def_Delta}}{>}0.
	\end{align*}
	Since $\barX(0) \in \Int \A_{\epsilon,0}$, and by \ref{hyp_C0_dA_nonlin}, we have a continuous map $t\mapsto d_{\partial\A_{\epsilon,t}}(\barX(t))$, $\barX(t) \in \Int\A_{\epsilon,t}$. So we may set $\b x_1(\cdot)\equiv \barX(\cdot)$ over $[t_0,t_1]$ and move to the next interval $[t_1,t_2]$.\\
	
	\noindent\textit{Case 2}: If $d_{\A_{\epsilon,t_0}}(\barX(t_0))\leq\eta/2$, let $\b u_0\in M_u\BB_M$ and $\b v_0=\b f(t_0,\barX(t_0),\b u_0)\in M_v\BB_N$ be as provided by \ref{hyp_IP_Keps_nonlin}. Let $\b y(\cdot)$ be defined on $[t_0,t_1]$ as follows
	\begin{equation}
	\b y(s):=\left\{\begin{array}{cl}
	\barX(t_0) + (s-t_0)\b v_0\; &\forall s\in[t_0,\min(t_0+k\barRho,t_1)],\\ 
	\barX(s-k\barRho)+k\barRho\b v_0\; &\forall s\in[\min(t_0+k\barRho,t_1),t_1].
	\end{array} \right. \label{eq_def_y}
	\end{equation}\\

	Applying classical measurable selection theorems for set-valued maps (see e.g.\ \cite{castaing_convex_1977}), fix a measurable selection $\b u(\cdot)$ satisfying \ref{hyp_absC_nonlin} as follows
	\begin{align*}
	&\bullet \text{[Case 2.1]}\; \forall\, s\in[t_0,\min(t_0+k\barRho,t_1)],\\
	&\|\b u(s)-\b u_0\|\le \beta_u(s),\, \\
	&\|\b f(s,\b y(s),\b u(s))-\b f(t_0,\b y(s),\b u_0)\|\leq \int_{t_0}^{s} \gamma\left(\sigma\right) d \sigma,\\
	&\bullet \text{[Case 2.2]}\; \forall\, s\in[\min(t_0+k\barRho,t_1),t_1],\\
	&\|\b u(s)-\barU(s-k\barRho)\|\le \beta_u(s),\, \\
	& \|\b f(s, \barX(s-k\barRho),\b u(s))-\b f(s-k\barRho,\barX(s-k\barRho), \barU(s-k\barRho))\|\\ 
	&\leq \int_{s-k\barRho}^{s} \gamma\left(\sigma\right) d \sigma,\\
	&\bullet \text{[Case 2.3]}\; \forall\, s\in[t_1,T],\; \b u(s) = \barU(s).
	\end{align*}
	Thus, at first we take the $t_0$-control as reference, then we change the reference to a delayed control, and finally we revert to the original control for the rest of the trajectory. Based on \ref{hyp_absC_nonlin}, we can bound the approximation error and prove that the control of \emph{Case 2.1} is applied long enough to make us stay in the constraint set for the whole $[t_0,t_1]$.

	Fix $\b x(t_0)=\barX(t_0)$ and use the control $\b u(\cdot)$ on $[t_0,T]$ to define a $\b f$-trajectory $\b x(\cdot)$. We have to show that $\b x(\cdot) \in \Int \A_{\epsilon,0}$ over $[t_0,t_1]$ and that it remains close to $\barX(\cdot)$ over $[t_0,T]$. We thus split the proof into three cases, one for each time interval in defining $\b u(\cdot)$.  To simplify our notations, we assume from now on that $t_0+k\barRho\leq t_1$ (otherwise Case 2.1 contains Case 2.2). We will proceed iteratively over $[t_i,t_{i+1}]$ so we must show that the assumptions on $\barX(\cdot)$ also hold for $\b x(\cdot)$.\\
	
	\begin{claim} We have the inclusion $\b x([t_0,T])\subset (R-1)\BB_N$ 	and for all $\delta\in[0,T]$ $$\sup_{\substack{[s,t]\subset[t_0,T]\\|s-t|\leq \delta}} \|\b x(t)-\b x(s)\|\leq \bar{\omega}(\delta).$$%
	\end{claim} %
	\noindent\textbf{Proof of Claim 2:} Note that by construction, since $\b u_0\in M_u\BB_M$, for all $s,t\in[0,T]$ with $s\le t$, by Cauchy-Schwarz (C-S) inequality, 
	\begin{align*}
		\int_{s}^{t} \theta(\sigma)\|\b u(\sigma)\|d\sigma&\le M_u\int_{s}^{t} \theta(\sigma)d\sigma +1_{t\ge t_1}\int_{\max(s,t_1)}^{t} \hspace{-1cm}\theta(\sigma) \|\barU(\sigma)\|d\sigma\\
		&\hspace{-1cm}+ 1_{s\le t_1} \int_{s}^{\min(t,t_1)} \hspace{-1cm}\theta(\sigma) (\|\barU(\sigma-k\barRho)+\beta_u(\sigma-k\barRho)\|)d\sigma\\
		&\hspace{-1cm}\le M_u\|\theta(\cdot)\|_{L^1(s,t)}\\
		&\hspace{-1cm}+\|\theta(\cdot)\|_{L^{2}(s,t)}(\|\barU(\cdot)\|_{L^{2}(0,T)}+\|\beta_u(\cdot)\|_{L^{2}(0,T)}).
	\end{align*}
	
	Hence, applying \ref{hyp_subLin_nonlin}, for any $t\in[t_0,T]$, with $s=t_0$,
	\begin{align*}
	\|\b x(t)\|&\leq \|\b x(t_0)\|+\int_{t_0}^{t} \|\b f(\sigma, \b x(\sigma),\b u(\sigma))\|d\sigma\\ &\leq\|\barX(\cdot)\|_{L^{\infty}(0,T)}+\int_{t_0}^{t} \theta(\sigma)(1+\|\b x(\sigma)\|+\|\b u(\sigma)\|)d\sigma\\
	\|\b x(t)\|&\le \|\barX(\cdot)\|_{L^{\infty}(0,T)}+(1+M_u)\|\theta(\cdot)\|_{L^{1}(0,T)}\\ &\hspace{0mm}+\|\theta(\cdot)\|_{L^{2}(t_0,t)}(\|\barU(\cdot)\|_{L^{2}(0,T)}+\|\beta_u(\cdot)\|_{L^{2}(0,T)})\\
	&+\int_{t_0}^{t} \theta(\sigma)\|\b x(\sigma)\|d\sigma
	\end{align*}
	Hence, owing to Gronwall's lemma, with $R$ as in \eqref{eq_def_R}, $\|\b x(t)\| \le R-1$. Moreover, for all $\delta>0$ and $[s,t]\subset[t_0,T]$ such that $|s-t|\leq \delta$
	\begin{align*}
	\|\b x(t)-\b x(s)\|&=\int_{s}^{t} \|\b f(\sigma, \b x(\sigma),\b u(\sigma))\|d\sigma\\ 
	& \le \int_{s}^{t} \theta(\sigma)(1+\|\b x(\sigma)\|+\|\b u(\sigma)\|)d\sigma\\
	\intertext{With the same computation as above}
	\|\b x(t)-\b x(s)\|& \le 	(1+R + M_u) \|\theta(\cdot)\|_{L^1(s,t)} \\ 
	&\hspace{-1.5cm}+\|\theta(\cdot)\|_{L^{2}(s,t)}(\|\barU(\cdot)\|_{L^{2}(0,T)}+\|\beta_u(\cdot)\|_{L^{2}(0,T)}) \le\bar{\omega}(\delta)
	\end{align*}
	\begin{flushright}
		$\square$
	\end{flushright}
	\noindent\textit{Case 2.1}: Let $t\in[t_0,t_0+k\barRho]$, by definition of $\b y(\cdot)$ in \eqref{eq_def_y},
	\begin{align}
	&\b x(t)- \b y(t)=\b x(t)- \b x(t_0)-(t-t_0)\b v_0\nonumber\\ 
	&\hspace{1cm}=\int_{t_0}^{t} \left[\b f(s, \b x(s),\b u(s))-\b f(t_0,\barX(t_0),\b u_0)\right]ds \label{eq_x-y_int}\\	
	&\|\b x(t)- \b y(t)\|\nonumber\\ 
	&\leq  \int_{t_0}^{t} \left\|\b f(s, \b x(s),\b u(s))-\b f(s,\b y(s),\b u(s))\right.\nonumber\\ 
	&\hspace{2cm}\left.+\b f(s,\b y(s),\b u(s))-\b f(t_0,\barX(t_0),\b u_0)\right\|ds\nonumber\\
	& \stackrel{\ref{hyp_lipF_nonlin}}{\le}	\int_{t_0}^{t} k_f(s)\|\b x(s)-\b y(s)\|ds\nonumber\\
	&\hspace{1cm} +\int_{t_0}^{t} \|\b f(s, \b y(s),\b u(s))-\b f(s,\barX(t_0),\b u(s))\|ds\nonumber\\
	&\hspace{1cm}+\int_{t_0}^{t}\|\b f(s,\barX(t_0),\b u(s))-\b f(s,\barX(t_0),\b u_0)\|ds.\nonumber\\
	\intertext{Apply now \ref{hyp_absC_nonlin}-\ref{hyp_lipF_nonlin} and note that $ [t_0,t_0+k\barRho]\subset [t_0,t_1]\subset [t_0,t_0+\Delta]$, so with $\omega_\gamma$ and $\omega_f$ as in \eqref{eq_def_omega}}
	&\|\b x(t)- \b y(t)\| \le (t-t_0) \omega_\gamma(\Delta)\\
	&\hspace{1cm}+\int_{t_0}^{t} k_f(s)\left[\|\b x(s)-\b y(s)\|+\|\b y(s)-\barX(t_0)\| \right]ds\nonumber\\
	&\le  (t-t_0)\omega_\gamma(\Delta)+\int_{t_0}^{t} k_f(s)\|\b x(s)-\b y(s)\|ds\nonumber\\
	&\hspace{1cm}+\int_{t_0}^{t} (t-t_0) M_v  k_f(s) ds\nonumber\\
	&\leq (t-t_0) (\omega_\gamma(\Delta) + M_v\omega_f(\Delta))\\
	&\hspace{1cm}+\int_{t_0}^{t} k_f(s)\|\b x(s)-\b y(s)\|ds.\nonumber
	\end{align}
	From Gronwall's lemma, with $M_{\Delta}$ and $C_{v,\Delta}$ as in \eqref{eq_def_csts_nonlin}, we deduce that $\|\b x(t)- \barX(t)\|$ is small:
	\begin{align*}
	\|\b x(t)- \b y(t)\|&\leq (t-t_0) M_{\Delta} e^{\int_{t_0}^{t}k_f(s)ds}\nonumber\\
	&\leq (t-t_0)M_{\Delta} e^{\omega_f(\Delta)} \le k\barRho C_{v,\Delta}\\
	\|\b x(t)- \barX(t_0)\|&\leq \|\b x(t)- \b y(t_0)\|+\|\b y(t)- \barX(t_0)\|\nonumber\\
	&\le (t-t_0)(M_v+M_{\Delta} e^{\omega_f(\Delta)})\\
	\|\b x(t)- \barX(t)\|&\leq \|\barX(t)- \barX(t_0)\|+\|\b x(t)- \barX(t_0)\| \nonumber\\
	&\le \bar{\omega}(k\barRho)+ k\barRho C_{v,\Delta}.
	\end{align*}
	Furthermore, setting  $$\varphi(t):=\int_{t_0}^{t} [\b f(s, \b x(s),\b u(s))-\b f(t_0,\barX(t_0),\b u_0)]ds,$$ 
	\begin{align*}	
	\|\b x(t)- \b x(t_0) - (t-t_0)\b v_0\|=\|\varphi(t)\|\nonumber\\
	&\hspace{-2cm}\leq (t-t_0)M_{\Delta} e^{\omega_f(\Delta)}\leq (t-t_0) \frac{\xi}{2}.
	\end{align*}
	In other words, $\b x(t)\in \b x(t_0)+(t-t_0)(\b v_0+\frac{\xi}{2}\BB_N)$ as in \eqref{eq_IPC_nonlin}. Hence, owing to \eqref{hyp_IP_Keps_nonlin},  since $t-t_0\le k\barRho\le \xi$, $\b x(t)\in \Int\A_{\epsilon, t}$.\\ 
	
	\noindent\textit{Case 2.2}: Let $t\in[t_0+k\barRho,t_1]$. As $\b x(t_0+k\barRho)= \b x(t_0) + k\barRho\,\b v_0 +\varphi(k\barRho)$, by \eqref{eq_def_y}, applying \eqref{eq_x-y_int},
	\begin{align*}
	\b x(t)- \b y(t)&=\b x(t) - \b x (t_0+k\barRho) +\barX (t_0)\nonumber\\
	&\hspace{1cm}+ \varphi(k\barRho)- \barX(t-k\barRho)\\
	\|\b x(t)- \b y(t)\|&\leq  \|\varphi(k\barRho)\|\\
	&\hspace{1cm}+\int_{t_0+k\barRho}^{t}\|\b x'(s)-\barX'(s-k\barRho)\|ds
	\end{align*}
	For any $s\in [t_0+k\barRho,t_1]$, by definition of $\b u(\cdot)$,
	\begin{align*}
	&\|\b x'(s)-\barX'(s-k\barRho)\|\le \int_{s-k\barRho}^{s} \gamma\left(\sigma\right) d \sigma\\
	&\hspace{1cm}+ \|\b f(s, \b x(s),\b u(s))-\b f(s,\barX(s-k\barRho), \barU(s-k\barRho))\|.
	\end{align*}
	However
	\begin{align*}
	&\int_{t_{0}+k \barRho}^{t} \int_{s-k \barRho}^{s} \gamma\left(\sigma\right) d \sigma d s \stackrel{\text{Fubini}}{=}\\
	&\int_{t_{0}}^{t} \gamma\left(\sigma\right)\int_{t_{0}}^{t_{1}} \chi_{\left[t_{0}+k \barRho, t\right]}(s) \chi_{[s-k \barRho, s]}\left(\sigma\right) d s  d \sigma \leq \omega_\gamma(\Delta) k \barRho,
	\end{align*}
	 hence by \ref{hyp_lipF_nonlin} and \eqref{eq_def_y}
	\begin{align*}
	&\int_{t_{0}+k\barRho}^{t}\|\b x'(s)-\barX'(s-k\barRho)\|ds\\
	&\leq \int_{t_{0}+k\barRho}^{t}k_f(s)(\|\b x(s)-\b y(s)\|+k\barRho M_v)ds+\omega_\gamma(\Delta) k \barRho
	\end{align*}\\

	Consequently, owing to Gronwall's lemma, and as $$\|\varphi(k\barRho)\|\leq k\barRho M_{\Delta} e^{\omega_f(\Delta)}$$ with the constants defined in \eqref{eq_def_csts_nonlin}
	\begin{align}
	\|\b x(t)- \b y(t)\|&\leq k\barRho (\omega_\gamma(\Delta)+M_v \omega_f(\Delta)\nonumber\\
	&\hspace{1cm}+ M_{\Delta} e^{\omega_f(\Delta)})e^{\omega_f(\Delta)}\nonumber\\
	&= k\barRho M_\Delta (1+e^{\omega_f(\Delta)}) e^{\omega_f(\Delta)} \label{eq_bound_x-y_2.2}
	\end{align}
	Hence, again by definition of $\b y(\cdot)$ in \eqref{eq_def_y},
	\begin{align*}
	\|\b x(t)- \barX(t)\|&\leq \|\b x(t)- \b y(t)\|+\|\b y(t)- \barX(t)\|\\
	&\leq \|\barX(t)- \barX(t-k\barRho)\|\\
	&\hspace{1cm}+  k\barRho [M_v +  M_\Delta (1+e^{\omega_f(\Delta)}) e^{\omega_f(\Delta)}]\\
	\|\b x(t)- \barX(t)&\|\leq \bar{\omega}(k\barRho)\\
	&\hspace{1cm}+ k\barRho (C_{v,\Delta}+ M_\Delta e^{2\omega_f(\Delta)}).
	\end{align*}
	Let $$\bm{\pi}(t):=\argmin_{\b z\in\A_{\epsilon,t-k\barRho}} \|\barX(t-k\barRho)-\b z\|.$$ Since $\barRho=\rho_{\epsilon,[0,T]}(\barX(\cdot))\le \hat{\rho}$ as in \eqref{eq_def_dN},
	\begin{align*}
	&\| \bm{\pi}(t)- \barX(t_0)\|\\
	&\leq \|\bm{\pi}(t)- \barX(t-k\barRho)\|+  \|\barX(t-k\barRho)- \barX(t_0)\|\nonumber\\
	&\leq \hat{\rho} + \bar{\omega}(\Delta)\leq \xi,
	\end{align*}\\

	Due to the inward-pointing condition \ref{hyp_IP_Keps_nonlin}, as $k\barRho\leq\xi$ and $\bm{\pi}(t)\in \A_{\epsilon,t-k\barRho}$, $$\bm{\pi}(t) + k\barRho(\b v_0 + \xi\BB_N)\subset \A_{\epsilon,t}.$$ Applying \eqref{eq_bound_x-y_2.2}, as by \eqref{eq_def_y} in this case $\b y(t)=\barX(t-k\barRho)+ k\barRho\b v_0$,
	\begin{align*}
	&\|\bm{\pi}(t) + k\barRho\b v_0 - \b x(t)\|\\
	&\leq \|\bm{\pi}(t)- \barX(t-k\barRho)\|+  \|\barX(t-k\barRho)+ k\barRho\b v_0 - \b x(t)\|\\
	&\leq \barRho + k\barRho M_\Delta (1+e^{\omega_f(\Delta)}) e^{\omega_f(\Delta)}\leq k \barRho\frac{\xi}{2},
	\end{align*}
	consequently $\b x(t)\in
	\bm{\pi}(t) + k\barRho(\b v_0 + \frac{\xi}{2}\BB_N)\subset \Int\A_{\epsilon,t}$.
	
	To summarize Cases 2.1 and 2.2, we have shown that for $t\in[t_0,t_1]$ we have $\b x(t)\in\Int\A_{\epsilon,t}$ and that
	\begin{equation}\label{eq_bound_2.1-2.2}
	\|\b x(t)- \barX(t)\|\leq 	\bar{\omega}(k\barRho)+ k\barRho (C_{v,\Delta}+ M_\Delta e^{2\omega_f(\Delta)}).
	\end{equation}
	
	\noindent\textit{Case 2.3}: Finally, let $t\in[t_1,T]$. As for $s\geq t_1$, $\b u(s) = \barU(s)$, by Filippov's theorem \cite{filippov_classical_1967} over $[t_1,T]$, with the bound \eqref{eq_bound_2.1-2.2} for $[t_0,t_1]$,
	\begin{align}
	&\forall\, \tau \in[t_0,T],\;\|\b x(\tau)- \barX(\tau)\|\nonumber\\
	&\leq e^{\omega_f(T)} \left[\bar{\omega}(k\barRho)+ k\barRho (C_{v,\Delta}+ M_\Delta e^{2\omega_f(\Delta)})\right]=g(\barRho). \label{eq_x_rec_infty}
	\end{align}\\

	\noindent\textbf{Step 2}: We iterate the procedure over $[t_0,t_{i}]$ for $i\in\iv{1}{N_0}$\\
	
	Set $\b x_0(\cdot):=\barX(\cdot)$, $\rho_0:=\barRho$, and $\b x_1(\cdot):=\b x(\cdot)$. Compute $\rho_1:=\rho_{\epsilon,[t_1,T]}(\b x_1(\cdot))$ and repeat the procedure described in Step 1 over  $\left[t_{1}, T\right]$ taking as reference trajectory $\b x_{1}(\cdot)$ to construct a $\b f$-trajectory $\b x_{2}(\cdot)$ strictly feasible over $\left[t_{1}, t_{2}\right]$. Set $\b x_{2}(\cdot)\equiv \b x_{1}(\cdot)$ on $[t_0,t_1]$. This iterative procedure defines a sequence of $\b f$-trajectories $\b x_i(\cdot)$. Each $\b x_i(\cdot)$ is strictly $\A_\epsilon$-feasible on $[t_0,t_{i}]$ and the corresponding $\b u_i(t)$ is always compared to the fixed $\barU(\cdot)$ on $[t_i,T]$ provided we can repeat the content of Step 1. Claim 1 holding also for $\b x_i(\cdot)$, we can indeed proceed with the iterations as long as $\rho_i:=\rho_{\epsilon,[t_{i},T]}(\b x_{i}(\cdot))$ is smaller than $\hat{\rho}$ and that $x_i([0,T])\subset 1+2 \|\barX(\cdot)\|_{L^{\infty}(0,T)}$ in order to apply \ref{hyp_absC_nonlin} to select the controls $\b u(\cdot)$. We prove these two claims by contradiction: fix $n$ as the first integer in $\iv{1}{N_0-1}$ such that $\rho_{n}>\hat{\rho}$ or that $x_n([0,T])\not\subset 1+2 \|\barX(\cdot)\|_{L^{\infty}(0,T)}$. Recall that $t_{N_0}=T$. We have for all $i\in \iv{0}{n-1}$
	\begin{align}
	\rho_{i+1}&\le \rho_{i} + \|\b x_{i+1}(\cdot)- \b x_{i}(\cdot)\|_{L^{\infty}(0,T)}\nonumber \\
	&\stackrel{\eqref{eq_x_rec_infty}}{\le}  \rho_i + g(\rho_i) = \tilde{g}(\rho_i).\label{eq_iter_rho}
	\end{align}
	Hence $d_i:=\|\b x_i(\cdot)- \barX(\cdot)\|_{L^{\infty}(0,T)}$ satisfies 
	\begin{align}
	d_{n}&\le\sum_{i=0}^{n-1} \|\b x_{i+1}(\cdot)- \b x_{i}(\cdot)\|_{L^{\infty}(0,T)} \stackrel{\eqref{eq_x_rec_infty}}{\le} \sum_{i=0}^{n-1}  g(\rho_i)\nonumber\\
	& \stackrel{\eqref{eq_iter_rho}}{\le}  \sum_{i=0}^{n-1}  \tilde{g}(\rho_i) \stackrel{\eqref{eq_iter_rho}}{\le} \sum_{i=0}^{n-1}  \tilde{g}^{\circ (i+1)}(\barRho)=\tilde{d}_{n}(\barRho),\label{eq_bound_dn_tdn}
	\end{align}
	where $\tilde{g}^{\circ i}$ corresponds to $\tilde{g}$ composed $i$-times with itself. Therefore we deduce that $$d_{n} \le \tilde{d}_{n}(\barRho)\leq \min(\lambda,\|\barX(\cdot)\|_{L^{\infty}(0,T)})$$ and, from \eqref{eq_def_dN}, that
	\begin{gather*}
	\rho_{n} \stackrel{\eqref{eq_iter_rho}}{\le} \tilde{g}^{\circ n}(\barRho)\leq \hat{\rho},\\ 
	\|\b x_n(\cdot)\|_{L^{\infty}(0,T)}\leq d_{n}+\|\barX(\cdot)\|_{L^{\infty}(0,T)}\leq 2 \|\barX(\cdot)\|_{L^{\infty}(0,T)}.
	\end{gather*}
	Consequently, $n$ as above does not exist and we can do the construction up to $\b x_{N_0}(\cdot)$. We have thus proved that 	
	$$ \|\b x_{N_0}(\cdot)- \barX(\cdot)\|_{L^{\infty}(0,T)}=d_{N_0}\leq \min(\lambda,\|\barX(\cdot)\|_{L^{\infty}(0,T)}).$$
	The trajectory $\b x^\epsilon(\cdot) \equiv \b x_{N_0}(\cdot)$ with control $\b u^\epsilon(\cdot)$ satisfies the requirements of the first statement of Theorem \ref{thm_NFT_unbounded}.\\
	
	To prove the second statement of Theorem \ref{thm_NFT_unbounded}, we consider additionally \ref{hyp_lipSelection_nonlin} to be satisfied.  Denote by $\|\cdot\|_{\b R(t)}$ the $\R^M$-norm $\|\b R(t)^{1/2}\cdot\|$. Since $\b R(\cdot)$ is continuous, it is uniformly continuous on the compact set $[0,T]$. Considering the operator norm, let $\bar{\mu}:=\|\b R(\cdot)\|^{1/2}_{L^\infty(0,T)}$ and $$\omega_{\b R}(\delta):=\sup_{[s,t]\subset[0,T], |s-t|\leq \delta} \|\b R(t)-\b R(s)\|.$$
	We now make further precise the approximation over the $L^2$-norm of the controls since, for all $i\in \iv{0}{N_0-1}$, if $d_{\partial\A_{\epsilon,t_i}}(\b x^\epsilon(t_i))\leq\eta/2$, then
	\begin{align*}
	&\forall\, s\in[t_i,t_i+k\rho_i],\; &&\|\b u^\epsilon(s)-\b u_i\|_{\b R(s)}\\
	& &&\hspace{-8mm}\leq (s-t_i)^\alpha\cdot\bar{\mu} k_u(s) \text{ for }  \|\b u_i\|\leq M_u,\\
	&\forall\, s\in[t_i+k\rho_i,t_{i+1}],\; &&\|\b u^\epsilon(s)-\barU(s-k\rho_i)\|_{\b R(s)}\\
	& &&\hspace{-8mm}\leq (k\rho_i)^\alpha\cdot\bar{\mu} k_u(s)\stackrel{k\rho_i\leq 1}{\leq} \bar{\mu}k_u(s);
	\end{align*}
	and, if $d_{\partial\A_{\epsilon,t_i}}(\b x^\epsilon(t_i))>\eta/2$, for all $s\in[t_i,t_{i+1}]$, $\b u(s) = \barU(s)$. Splitting the intervals of integration,
	\begin{align*}
	&\int_0^T \Vert\barU(t)\Vert_{\b R(t)}^2-\Vert\b u^\epsilon(t)\Vert_{\b R(t)}^2 dt\\
	&= \int_0^T \Vert\barU(t)^2\Vert_{\b R(t)} dt\\
	&+\sum_{i=0}^{N_0-1} \int_{t_i}^{t_i+k\rho_i}[\Vert\b u_i\Vert_{\b R(t)}^2-\Vert\b u^\epsilon(t)\Vert_{\b R(t)}^2] dt\\
	&+\sum_{i=0}^{N_0-1} \int_{t_i+k\rho_i}^{t_{i+1}} [\Vert\barU(t-k\rho_i)\Vert_{\b R(t)}^2-\Vert\b u^\epsilon(t)\Vert_{\b R(t)}^2] dt\\
	&-\sum_{i=0}^{N_0-1} \left[\int_{t_i}^{t_i+k\rho_i} \Vert\b u_i\Vert_{\b R(t)}^2 dt+ \int_{t_i+k\rho_i}^{t_{i+1}} \Vert\barU(t-k\rho_i)\Vert_{\b R(t)}^2 dt\right].
	\end{align*}
	Since 
	\begin{align*}
	\Vert\barU(t-k\rho_i)\Vert_{\b R(t)}^2&= \barU(t-k\rho_i)^\top \left(\b R(t-k\rho_i)+ \b R(t)\right.\\
	&\hspace{5mm}\left.- \b R(t-k\rho_i)\right) \barU(t-k\rho_i),
	\end{align*}
	the last term cancels most of the first term. Using that, for any $a,b\in\R_+$, $|a^2-b^2|\leq |a-b|\,(2b+|a-b|)$, we can also bound the integrals containing differences,
	\begin{align*}
	&\left\vert\int_0^T \Vert\barU(t)\Vert_{\b R(t)}^2-\Vert\b u^\epsilon(t)\Vert_{\b R(t)}^2 dt\right\vert\\
	&\leq \sum_{i=0}^{N_0-1} \int_{t_{i+1}-k\rho_i}^{t_{i+1}} \Vert\barU(t)\Vert_{\b R(t)}^2 dt\\
	&+ \sum_{i=0}^{N_0-1} \int_{t_i}^{t_i+k\rho_i}\bar{\mu}^2(k\rho_i)^\alpha k_u(t)(2M_u+k_u(t)) dt\\
	&+\sum_{i=0}^{N_0-1} \int_{t_i+k\rho_i}^{t_{i+1}} \bar{\mu}^2(k\rho_i)^\alpha k_u(t)(2\Vert\barU(t-k\rho_i)\Vert+k_u(t)) dt\\
	&+\sum_{i=0}^{N_0-1} \int_{t_i+k\rho_i}^{t_{i+1}} \Vert\b R(t)- \b R(t-k\rho_i)\Vert\cdot\Vert\barU(t-k\rho_i)\Vert^2 dt \\
	&+\sum_{i=0}^{N_0-1} k\rho_i M_u^2 \bar{\mu}^2.
	\end{align*}
	As $\rho_i\le \tilde{g}^{\circ N_0}(\barRho)$ and,by \eqref{eq_bound_dn_tdn}, $\sum_{i=0}^{N_0-1} \rho_i\le \tilde{d}_{N_0}(\barRho)$, we can bound all the terms $k\rho_i$ in order to derive a r.h.s.\ only depending on $\barRho$,
	\begin{align*}
	&\left\vert\int_0^T \Vert\barU(t)\Vert_{\b R(t)}^2-\Vert\b u^\epsilon(t)\Vert_{\b R(t)}^2 dt\right\vert\\
	&\leq \sum_{i=0}^{N_0-1} \int_{t_{i+1}-\tilde{g}^{\circ N_0}(\barRho)}^{t_{i+1}} \bar{\mu}^2\Vert\barU(t)\Vert^2 dt\\
	& + \bar{\mu}^2k^\alpha\,\tilde{g}^{\circ N_0}(\barRho)^\alpha\left[2M_u\|k_u(\cdot)\|_{L^1(0,T)} +\|k_u(\cdot)\|_{L^2(0,T)}^2\right] \\
	& + 2\bar{\mu}^2k^\alpha\,\tilde{g}^{\circ N_0}(\barRho)^\alpha\|k_u(\cdot)\|_{L^2(0,T)}\|\barU(\cdot)\|_{L^2(0,T)}\\
	&+\omega_{\b R}(\tilde{g}^{\circ N_0}(\barRho)) \|\barU(\cdot)\|_{L^2(0,T)}^2+k\,\tilde{d}_{N_0}(\barRho) M_u^2\bar{\mu}^2.
	\end{align*}
	Taking limit in the r.h.s.\ when $\epsilon$ goes to $0^+$, we deduce that we can take $\epsArr$ small enough as to have a small $\barRho$ and the r.h.s.\ smaller than $\lambda$, which concludes the proof.  
	\begin{flushright}
		$\blacksquare$
	\end{flushright}	%

	\bibliography{IFAC_CAO_22} 

\end{document}